\documentclass[lettersize,journal]{IEEEtran}
\usepackage{amsmath,amsfonts}
\usepackage{algorithmic}
\usepackage{algorithm}
\usepackage{array}
\usepackage{textcomp}
\usepackage{stfloats}
\usepackage{url}
\usepackage{verbatim}
\usepackage{graphicx}
\usepackage{cite}

\usepackage{pifont}
\usepackage{amsmath,amssymb,amsfonts}
\usepackage{enumerate}
\usepackage{subfigure}
\usepackage{times}
\usepackage{float}
\usepackage{color}
\usepackage{cuted}
\begin{document}

\title{Temporal-Assisted Beamforming and Trajectory Prediction in Sensing-Enabled UAV Communications }

\author{\IEEEauthorblockN{
            Shengcai Zhou, \emph{Student~Member, IEEE},
            Halvin Yang,
            \emph{Member, IEEE},
			Luping Xiang, \emph{Member, IEEE},
			and~Kun~Yang}, \emph{Fellow, IEEE}
			\vspace{-0.5 cm}\\

       \thanks{The authors would like to thank the financial support of Quzhou Government (Grant No.: 2023D005), Natural Science Foundation of China (Grant No. 62132004 and Grant No. 62301122), and the Xiaomi Young Scholars - Technology Innovation Award. \textit{(Corresponding author: Luping Xiang.)}}
        \thanks{Shengcai Zhou is with the School of Information and Communication Engineering, University of Electronic Science and
        Technology of China, Chengdu 611731, China, email: 202222011112@std.uestc.edu.cn.}
        \thanks{Halvin Yang is with  the Department of Electronic and Electrical Engineering, University College London, WClE 7JE London, U.K, e-mail:uceehhy@ucl.ac.uk.}
        \thanks{Luping Xiang and Kun Yang are with the State Key Laboratory ofNovel Software Technology, Nanjing University, Nanjing 210008, China, and School of Intelligent Software and Engineering, Nanjing University (Suzhou Campus), Suzhou, 215163, China, email: luping.xiang@nju.edu.cn, kunyang@nju.edu.cn.}
	}

\maketitle



\maketitle

\begin{abstract}
In the evolving landscape of high-speed communication, the shift from traditional pilot-based methods to a Sensing-Oriented Approach (SOA) is anticipated to gain momentum. This paper delves into the development of an innovative Integrated Sensing and Communication (ISAC) framework, specifically tailored for beamforming and trajectory prediction processes. Central to this research is the exploration of an Unmanned Aerial Vehicle (UAV)-enabled communication system, which seamlessly integrates ISAC technology. This integration underscores the synergistic interplay between sensing and communication capabilities. The proposed system initially deploys omnidirectional beams for the sensing-focused phase, subsequently transitioning to directional beams for precise object tracking. This process incorporates an Extended Kalman Filtering (EKF) methodology for the accurate estimation and prediction of object states. A novel frame structure is introduced, employing historical sensing data to optimize beamforming in real-time for subsequent time slots, a strategy we refer to as 'temporal-assisted' beamforming. To refine the temporal-assisted beamforming technique, we employ Successive Convex Approximation (SCA) in tandem with Iterative Rank Minimization (IRM), yielding high-quality suboptimal solutions. Comparative analysis with conventional pilot-based systems reveals that our approach yields a substantial improvement of 156\% in multi-object scenarios and 136\% in single-object scenarios.
\end{abstract}


\section{Introduction}

\subsection{Background}
\IEEEPARstart{T}{he} evolution of digital vehicular technology has underscored the critical role of vehicle-to-everything (V2X) communications in societal frameworks. Accelerated by advances in 5G networks, the demand for V2X capabilities has increased, particularly for features such as low-latency transmissions, precise location tracking during high-speed movements, and robust anti-jamming mechanisms, which have traditionally been associated with substantial hardware investments \cite{8246850}. The recent discourse around the coexistence of radar and communication frequencies \cite{9173030,8288677} has paved the way for integrated sensing and communication (ISAC) technologies, carving out a novel research avenue \cite{8828023,9606831,8999605,9945983,10693589}. ISAC stands out by offering concurrent sensing and communication functionalities, thereby enhancing spectrum and energy efficiency. With integration into millimeter-wave and massive MIMO systems, ISAC signals are poised to deliver elevated data transmission rates and beamforming capabilities. Reciprocal improvement of sensing and communication operations facilitates synergistic improvement in performance, generating collaborative advantages \cite{9737357}. By utilising information derived from sensing, conventional pilot-based beam training can be substituted, significantly reducing the overhead of the system \cite{9171304}. Consequently, sensing-assisted communication has gained substantial traction within the V2X domain.

Conventional V2X scenarios, predominantly terrestrial, encounter limitations such as obstructions and insufficient coverage from the base station, which can severely hinder sensing operations. With the maturation of the wireless communication technology of unmanned aerial vehicles (UAVs) \cite{9456851,8709739}, UAVs have emerged as promising aerial platforms that are poised to overcome these challenges. Their inherent agility and mobility \cite{8870206,9151975} enable UAVs to access areas with poor connectivity and bypass numerous ground obstacles to establish direct line-of-sight (LoS) links \cite{9916163,9696263}. This is particularly advantageous in post-disaster scenarios, where UAVs can provide critical communication links and conduct search and rescue operations.  Therefore, the UAV-based ISAC system is envisioned as a significant adjunct to terrestrial ISAC infrastructures.

\subsection{Related Work}
Recent advancements in V2X communications have increasingly addressed the challenges posed by high-speed environments \cite{9201355}. High mobility contexts require frequent updates in channel estimation and beam alignment, traditionally relying on pilot signals, which incurs considerable overhead \cite{8851151}. Although the RAF algorithm has been proposed to mitigate pilot overhead while preserving data bandwidth and error probability, the inherent trade-off between the timeliness and accuracy of channel estimation against pilot costs persists. The introduction of orthogonal time frequency space (OTFS) modulation \cite{7925924} has improved the conventional limitations of OFDM, improving the effectiveness of the pilot estimation. ISAC-assisted OTFS transmission methods \cite{9557830} have built on this advancement, although initial estimation protocols remain unspecified. Alternative approaches, such as those presented in \cite{9171304}, leverage deep sensing integration with communication for efficient beam tracking, using radar attributes to replace pilot signals and uplink feedback, and introducing extended Kalman filter techniques for beam prediction, thus reducing communication latency. In addition to this, a novel message-passing algorithm is proposed in \cite{9246715} to reduce computational demands. Additionally, \cite{10330059} employs direction-of-arrival (DOA) information through reconfigurable intelligent sensing surfaces (RISS) instead of pilot signals. Addressing the challenges of directing narrow massive MIMO beams, studies such as \cite{9747255} and \cite{9947033} explore target enlargement strategies, including beam width modulation and stage-specific beamwidth selections. The ground-level complexities are further addressed by sensing-assisted beamforming techniques in \cite{9833957} and \cite{10061429}, which enhance geometric modeling research.  \textbf{However, the interrelation between omnidirectional and directional waveforms in the context of ISAC waveform design has not been adequately addressed.}

In UAV research that integrates ISAC, motion dynamics is a key focus. The stationary and mobility scenarios of UAVs are contemplated in \cite{9916163}, examining the co-optimisation of the communication and flight paths of UAVs within sensory constraints. Based on this, \cite{9739676} suggests strategies to equilibrate the durations of sensing and communication, thus improving sensory adaptability without excessive energy expenditure. Beyond mere beam gain metrics\cite{10313997}, \cite{10529184} investigates the trade-offs between minimising Cramer-Rao bound (CRB) and trajectory optimization, achieving superior sensing accuracy at a marginal communication performance cost. Additionally, the sequential implementation of projects involving multiple UAVs is discussed. \cite{9293257} addresses the challenges in UAV joint localization, user association, and power control for UAV transmissions. However, the exploration of multi-array beams and dynamic targets as noted in \cite{10529184,9293257} remains an area of improvement. Concerning ISAC framework protocols, \cite{10098686} delineates three prevalent designs for UAVs: Co-ISAC, TDM-ISAC, and Hybrid-ISAC, discussing their performance variations based on factors like service quality requirements, target positioning, and mobility. \textbf{The enhancement of ISAC performance remains an area for further investigation.} Recent work in \cite{10404096} integrates OTFS with UAV systems to propose novel power distribution strategies for UAV-user communications. \textbf{However, using temporal connections to further optimise beamforming matrix is not taken into account}. Our work is therefore inspired by the promising synergy between ISAC base stations and ISAC-equipped UAVs.

\begin{table*}[t]
\centering
\caption{Contrasting Our Contributions To The State-Of-Art}
\setlength{\tabcolsep}{4mm}{
\begin{tabular}{l|c|c|c|c|c|c|c|c}
 \hline 
 Contributions & \textbf{this work} & \cite{8288677} & \cite{9171304} & \cite{9916163}&\cite{10529184} & \cite{8851151} &\cite{9739676}&\cite{10404096} \\ \hline\hline
 ISAC & \ding{52} & \ding{51} & \ding{51} & \ding{51} &\ding{51} &&\ding{51}&\ding{51} \\
 \hline
 Pilot free & \ding{52} & \ding{51} &  & \ding{51}&\ding{51}&&\ding{51}&\ding{51}  \\
 \hline
 Beam coverage & \ding{52} & \ding{51} &  & \ding{51} &  &&& \\
 \hline
 Temporal-assisted & \ding{52} &  & \ding{51} & & \ding{51} &\ding{51}&& \\
 \hline
 Beam prediction & \ding{52} &  & \ding{51} &  & &\ding{51}&&  \\
 \hline
 Moving target & \ding{52} &  & \ding{51} &  & &\ding{51}&&\ding{51}  \\
 \hline
\end{tabular}
}
\vspace{-0.3 cm}
\label{contributions}
\end{table*}

\subsection{Our Contributions}
In this paper, we explore an emergency scenario where a UAV is deployed to establish a temporary communication network following the failure of a conventional system. The UAV moves continuously, seeking potential objects to establish communication links. Specifically, we introduce a novel ISAC framework mounted on UAVs, designed for optimized beam transmission and object tracking. This architecture innovatively employs an omnidirectional integrated beam for the initial sensing phase, eschewing conventional beam training and pilot reliance. Following the initial radar echo capture, we refine the process by deploying directional beams for enhanced object tracking, thereby securing more precise echo data, which is iteratively refined throughout the frame's duration. We refer to this process as temporal assistance. Subsequent frames involve re-initiating this omnidirectional emission to detect any new objects entering the monitored zone. To augment accuracy and minimize latency, we advance an extended Kalman filtering (EKF) methodology for meticulous tracking and prediction of object motion. The robustness of EKF lends superior anti-interference properties to beam alignment, maintaining high accuracy even amidst sporadic communication disruptions, courtesy of its multi-step predictive capability. In general, the proposed method not only includes the assistance of sensing to communication, but also includes the assistance of communication to sensing, so that the sensing accuracy and the communication rate can be improved simultaneously. The novelty of this work is explicitly contrasted in Table \ref{contributions}, and our contributions are encapsulated as follows:

\begin{itemize}
\item We present a UAV-centric, sensing-aided communication strategy within the ISAC system. A novel frame structure is conceptualized where the UAV initially deploys an omnidirectional ISAC beam to probe the environment, then harnesses the gathered data on objects and channels to craft a directional beam, thereby enhancing communication throughput and accruing precise information for subsequent optimization until the frame concludes.

\item An EKF-based approach is articulated for object state estimation and prediction, seamlessly integrating into the temporal-assisted beamforming workflow. We employ the method of successive convex approximation (SCA) \cite{Dinh2010LocalCO} alongside the iterative rank minimization (IRM) technique \cite{Sun2017TwoAF} to derive solutions of elevated quality.

\item Through simulation analyses, our proposed methodology exhibits superior performance compared to conventional pilot-based approaches. In a single frame, our approach demonstrates a notable performance advantage over the pilot-based method in time slot 6, achieving an improvement of $156\%$ for scenarios involving multiple objects and $136\%$ for a scenario with a single object. Furthermore, we scrutinize the comparative efficacy of our scheme against the Water-Filling algorithm across varying distances, underscoring our scheme's robustness in the face of significant parameter estimation discrepancies.
\end{itemize}
\begin{figure}[h] 
\centerline{\includegraphics[width=0.45\textwidth]{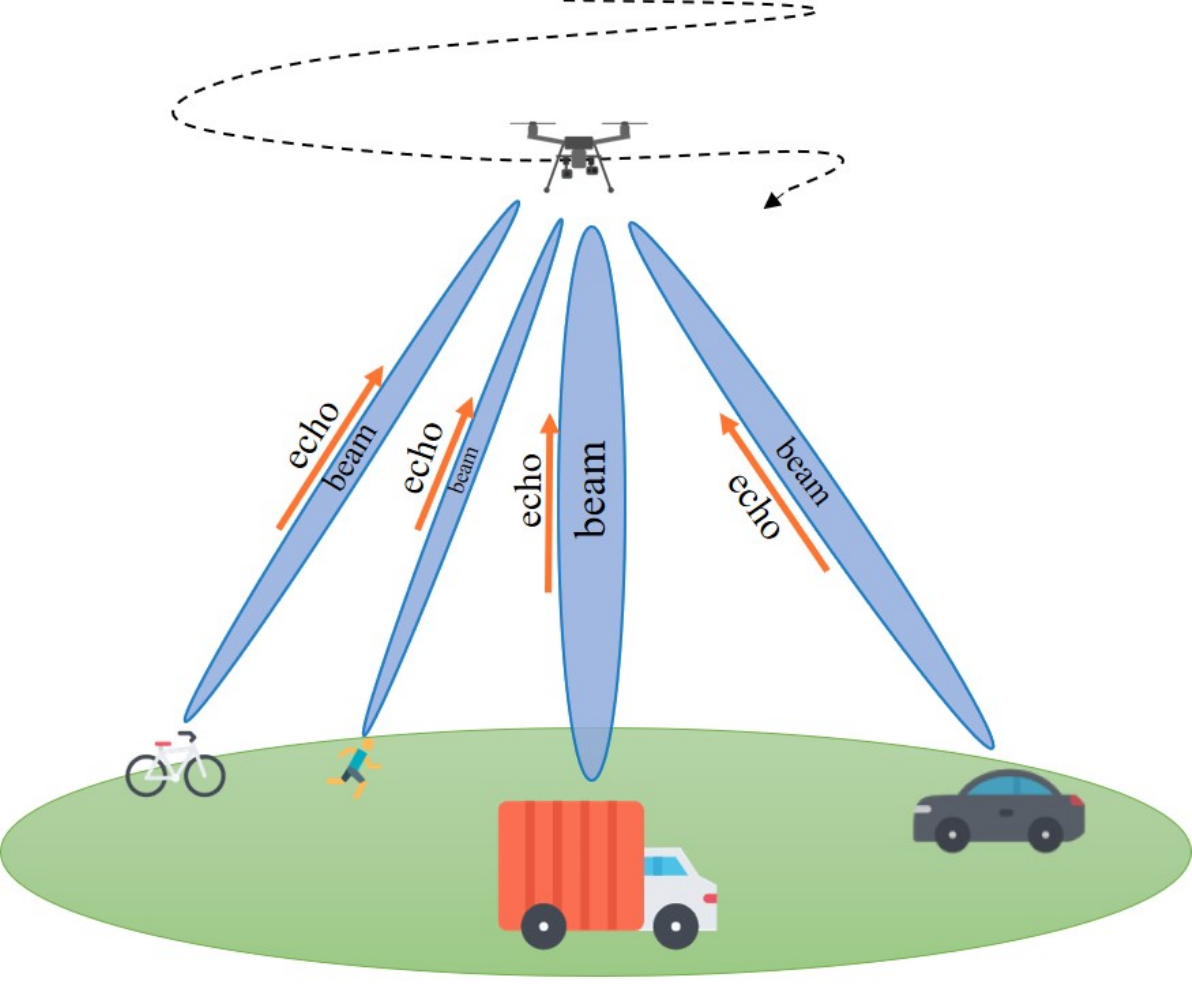}}
\caption{Sensing-enabled UAV system model.}
\label{fig.scenario}
\vspace{-0.3 cm}
\end{figure}

The structure of the remainder of this paper is as follows: Section II delineates the system model; Section III elaborates on the formulation of the optimization problem; the proposed solutions are expounded in Section IV; Section V is dedicated to the presentation of numerical results; and a conclusion is drawn in Section VI.

\section{System Model}

\subsection{Overall Structure}
As shown in Fig. \ref{fig.scenario}, we consider a UAV-enabled ISAC system, which serves $K$ objects on the road. To simplify the representation without losing generality, We assume that the UAV is equipped with a transmitting uniform linear array (ULA) containing $N_t$ antennas and a different receive ULA composed of $N_r$ antennas. Since the LoS component is 20 dB higher than the non-LoS component \cite{8851151}, it is hypothesized that the UAV engages in LoS communication with each object, which is conceptually considered a point target. The LoS MIMO channel model has few sub-channels, so the maximum delay spread is low with inconspicuous multipath effect.

In order to establish a reliable communication link, it is imperative for the UAV to ascertain precise data regarding the object's azimuth, range, and velocity, necessitating beam alignment towards the object to the highest degree feasible. conventional beam training practices are hampered by substantial delays and elevated communication overhead. In pursuit of heightened sensing accuracy and reduced overhead, an EKF predictive framework is integrated, and a novel sensing-assisted communication approach is introduced for the UAV-centric ISAC system.

The spatiotemporal parameters of the $k$-th object, namely the angle $\theta_k(t)$, distance $d_k(t)$, and velocity $v_k(t)$, are treated as time-dependent functions within the time $t\in[0,T]$, where $T$ is the maximum duration of the whole frame. For convenience, the time period $T$ is discretized into slots $\Delta T$, with $\theta_{k,n}$, $d_{k,n}$ and $v_{k,n}$ defining the motion metrics for each object in the $n$-th slot. Consistent with established theoretical models \cite{8025577}, it is assumed these parameters remain invariant within each $\Delta T$ slot.

\begin{figure}[h] 
\centerline{\includegraphics[width=0.45\textwidth]{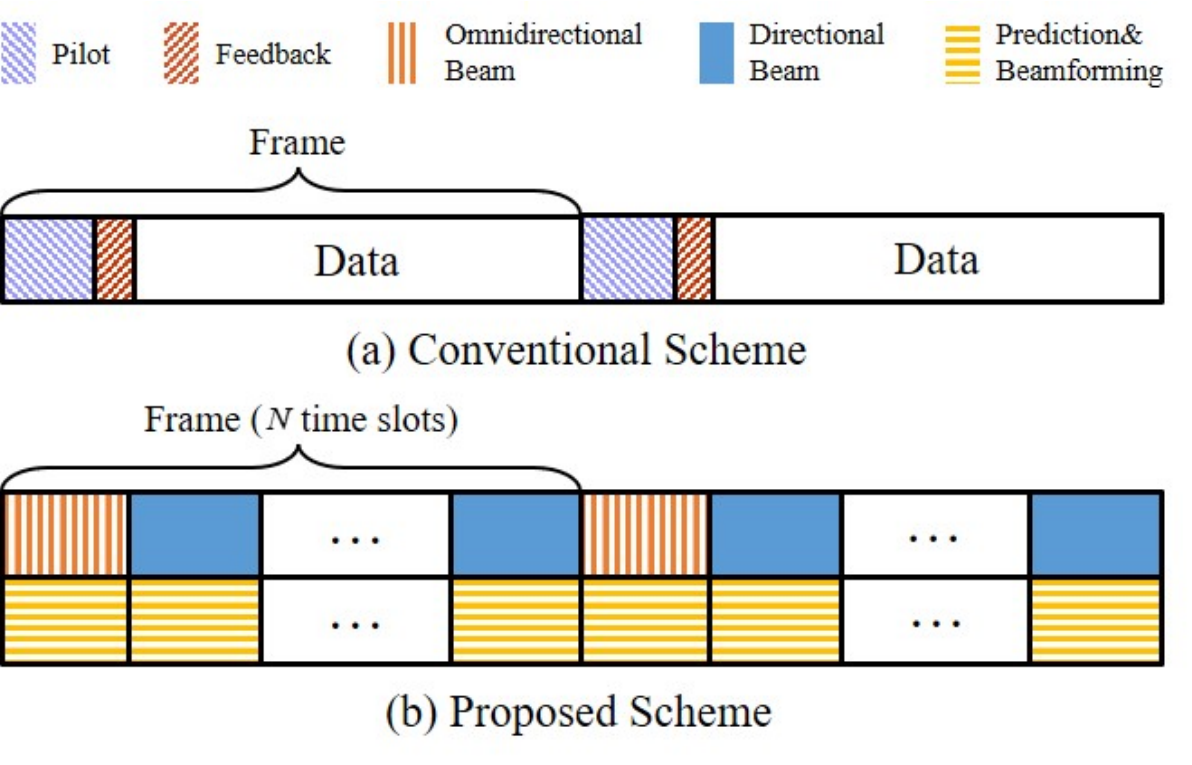}}
\caption{Frame structure.}
\label{fig.frame}
\vspace{-0.3 cm}
\end{figure}
\subsubsection{Frame Structure}
A comparative overview of the frame structures for conventional beam training and the proposed framework is delineated in Fig. \ref{fig.frame}. conventional methods involve the UAV transmitting a pilot signal to capture object parameters, followed by feedback for beamforming and data transmission. Conversely, the proposed framework utilizes the entire frame for simultaneous sensing and communication, obviating the need for dedicated downlink pilots and supplanting uplink feedback with radar echoes, thus enhancing communication efficacy while curtailing temporal overhead.

Specifically, frame interconnections are facilitated by periodic omnidirectional beams\footnote{In wide-area space, omnidirectional beams are transmited for sensing unkown objects, with carrying information to $K$ known objects. The channel information of all objects is then collected and used to optimize the communication rate. If $K=0$, the omnidirectional beams only execute sensing.}. These beams capture data on both newly sensed objects and those previously known. For new objects, directional beams are crafted based on the sensed data and refined over time. For previously known objects, having already secured sufficient data, we execute a two-step prediction based on the object's last known position from the previous frame, thereby swiftly restoring communication rates to near-optimal levels. In scenarios with high object maneuverability, the proposed beam prediction exhibits superior stability.

\subsubsection{Initial Estimation}
Diverging from conventional beam training, the initial estimation phase employs omnidirectional beams to detect objects, akin to a MIMO radar with communicative capabilities, with subsequent frame updates for new object detection.
\subsubsection{Directional Estimation}
During the $n$-th interval, the UAV assesses measured values $\hat{\theta}_{k,n}$, $\hat{\tau}_{k,n}$ and $\hat{\mu}_{k,n}$,  integrating them with predicted values, $\bar{\theta}_{k,n\vert n-1}$, $\bar{d}_{k,n\vert n-1}$ and $\bar{v}_{k,n\vert n-1}$, into the EKF for refined estimates $\bar{\theta}_{k,n}$, $\bar{d}_{k,n}$ and $\bar{v}_{k,n}$. Post-optimization, the UAV proceeds with beamforming, followed by the reacquisition of echoes, perpetuating this cycle. This prediction process is governed by kinematic equations, with the object's motion parameter evolution detailed in Section II-D.

\subsection{Communication Model}
Considering the UAV's awareness of $K$ objects within its operational range, it commences by formulating an ISAC signal $\mathbf{s}(t)$ in a $K$-dimensional space, which is defined as follows
\begin{equation}
\mathbf{s}_n(t)=[s_{1,n}(t),...,s_{K,n}(t)]^T,
\end{equation}
where the $k$-th signal $s_{k,n}(t)$ includes the communication information of the $k$-th object at time slot $n$, with independent circularly symmetric complex Gaussian (CSCG) random variables with zero mean and unit variance, represented by $s_{k,n}(t)\sim\mathcal{CN}(0,1)$.

Subsequently, after beamforming, the transmitted signal becomes
\begin{equation}
\mathbf{x}_n(t)=\mathbf{W}_n\mathbf{s}_n(t),
\end{equation}
where $\mathbf{W}_n\in \mathbb{C}^{N_t\times K}$ denotes the beamforming matrix, also referred to as the linear precoder.
\subsubsection{Omnidirectional Phase}
According to literature \cite{9916163}, the channel vector from the UAV to the $k$-th object is represented as
\begin{equation}
\mathbf{h}_{k,n}=\frac{\alpha_0}{d_{k,n}}e^{j\frac{2\pi}{\lambda}d_{k,n}}\mathbf{a}(\theta_{k,n}),
\end{equation}
where the terms $\alpha_0$ and $\lambda$ signify the path loss at a reference distance $d=1$m and carrier wavelength. The distance between the UAV and the $k$-th object at time slot $n$ is $d_{k,n}$. The steering vector aimed at object $k$ is expressed as
\begin{equation}
\mathbf{a}(\theta_{k,n})=[1,e^{j2{\pi}\frac{\Delta d}{\lambda}\sin{\theta_{k,n}}},...,e^{j2{\pi}\frac{\Delta d}{\lambda}(N_t-1)\sin{\theta_{k,n}}} ]^T,
\end{equation}
where $\Delta d$ represents the gap between two neighboring antennas and $\theta_{k,n}\in(-\pi,\pi)$ is the $k$-th object's angular position.

Consequently, the signal received by the $k$-th object is
\begin{equation}
y_{k,n}(t)=\mathbf{h}_{k,n}^He^{j2\pi\omega_{k,n}t}\mathbf{x}_n(t)+z_C(t),
\end{equation}
with $z_C(t)\sim\mathcal{CN}(0,\sigma_C^2)$ symbolizing the additive white Gaussian noise (AWGN) at the receiver. $\omega_{k,n}$ represents the Doppler frequency.

Based on the principles of omnidirectional beampattern \cite{4516997} with a unit total transmit power of $P_T$, the orthogonal nature of the signal $\mathbf{x}_n(t)$ is ensured. For scenarios where $K\geq N_t$, the spatial covariance matrix of the omnidirectional beam is
\begin{equation}
\mathbf{C}=\mathbb{E}[\mathbf{x}_n(t)\mathbf{x}^H_n(t)]=\mathbf{W}_n\mathbf{W}^H_n=\frac{P_T}{N_t}\mathbf{I}_{N_t},
\end{equation}
where $\mathbf{I}_{N_t}$ denotes the $N_t\times N_t$ identity matrix, the received SINR of the $k$-th object is represented as
\begin{equation}\label{omnisinr}
\gamma_{k,n}=\frac{\frac{1}{K}\mathbf{h}_{k,n}^H \mathbf{C} \mathbf{h}_{k,n}}{\frac{K-1}{K}\mathbf{h}_{k,n}^H \mathbf{C} \mathbf{h}_{k,n}+\sigma_C^2}=\frac{\frac{\alpha_0^2P_T}{d_{k,n}^2K}}{\frac{\alpha_0^2P_T(K-1)}{d_{k,n}^2K}+\sigma_C^2}.
\end{equation}

For cases where $K < N_t$, the linear precoder is deemed non-essential, allowing for the direct engagement of $K$ antennas to individually broadcast the ISAC signals. The SINR in this instance mirrors that of (\ref{omnisinr}).

Summarily, the attainable sum rate for the object $k$ during the $n$-th time slot is represented as
\begin{equation}\label{sumrate}
R_{k,n}=\log_{2}{(1+\gamma_{k,n})}.
\end{equation}

\subsubsection{Directional Phase}

The received signal for object $k$ is represented as
\begin{equation}
y_{k,n}(t)=\mathbf{h}_{k,n}^He^{j2\pi\omega_{k,n}t}\mathbf{W}_n\mathbf{s}_n(t)+z_C(t).
\end{equation}
In concise terms, the SINR for object $k$ can be articulated as
\begin{align}\label{dirsinr}
\gamma_{k,n}=\frac{|\mathbf{h}_{k,n}^H \mathbf{w}_{k,n}|^2}{\sum_{i=1,i\neq k}^K |\mathbf{h}_{k,n}^H \mathbf{w}_{i,n}|^2+\sigma_C^2},
\end{align}
where the beamforming vector, denoted as $\mathbf{w}_{k,n}$, stands as the $k$-th column vector extracted from $\mathbf{W}_n$. For time slot $n$, the achievable sum rate pertaining to  object $k$ aligns with expression $(\ref{sumrate})$.

\subsection{Radar Model}
\subsubsection{Omnidirectional Phase}
Omnidirectional radar is utilized for initial detection. Its orthogonal waveforms, when transmitted through antennas, ensure that the reflected signals from different objects are linearly independent, enabling the use of Capon and other adaptive array algorithms to distinguish between objects  \cite{4350230, 7071108}. For the $k$-th object, the radar echo is mathematically represented as
\begin{equation}
\mathbf{r}_{k,n}(t)=\beta_{k,n}e^{j2\pi\mu_{k,n}t}\mathbf{b}(\theta_{k,n})\mathbf{a}^H(\theta_{k,n})\mathbf{x}_n(t-\tau_{k,n})+\mathbf{z}_{k,n}(t),
\end{equation}
where $\mathbf{z}_{k,n}(t)$ represents the zero-mean complex additive white Gaussian noise with variance $\sigma^2$. The terms $\beta_{k,n}$ and $\tau_{k,n}$ denote the reflection coefficient and time delay for the $k$-th object during the $n$-th time slot, respectively, and $\mu_{k,n}$ is indicative of the Doppler frequency. The reflection coefficient is expressed as
\begin{equation}
\beta_{k,n}=\frac{\varepsilon}{d_{k,n}^2},
\end{equation}
in which $\varepsilon$ represents the complex radar cross-section (RCS). For analytical ease, we assume a constant RCS, typifying a Swerling I target\cite{Richards2005FundamentalsOR}.

Employing a standard matched-filtering approach \cite{6324717}, the system can estimate signal delays and Doppler frequencies. The resultant processed vectors are
\begin{equation}
\Tilde{\mathbf{R}}_{k,n}=G_m\beta_{k,n}\frac{P_T}{N_t}\mathbf{b}(\theta_{k,n})\mathbf{a}^H(\theta_{k,n})+\mathbf{Z}_r,
\end{equation}

where $G_m$ signifies the signal processing gain from matched-filtering. The noise matrix $\mathbf{Z}_r$ is composed of independent, zero-mean, complex Gaussian elements with variance $\sigma_r^2$. The radar echo SNR is calculated as
\begin{equation}
SNR_{k,n}=
\frac{ P_T G_m |\beta_{k,n}|^2}{\sigma^2}.
\end{equation}
\subsubsection{Directional Phase}
As per \cite{9557830}, steering vectors from different angles are asymptotically orthogonal under the massive MIMO regime, which generally prevents interference among reflected echoes. In instances where targets are closely spaced and this theory is inadequate, digital beamforming techniques can be applied to reduce echo interference. Similarly, the received echo for the directional phase is formulated as
\begin{align}
&\mathbf{r}_{k,n}(t) \nonumber \\
& =\beta_{k,n}e^{j2\pi\mu_{k,n}t}\mathbf{b}(\theta_{k,n})\mathbf{a}^H(\theta_{k,n})\mathbf{W}_n \mathbf{s}(t-\tau_{k,n}) + \mathbf{z}_{k,n}(t).
\end{align}
Finally, the SNR for the radar's echoed signal in this phase is described by
\begin{align}
SNR_{k,n}=\frac{G_m |\beta_{k,n}|^2 \sum_{i=1}^K|\mathbf{a}^H(\theta_{k,n})\mathbf{w}_{i,n}|^2 }{\sigma^2}.
\end{align}

\subsection{Measurement Model}

Upon the application of matched-filtering, the measurement model for ascertaining the distance $d_{k,n}$ and Doppler shift $\mu_{k,n}$ takes the following form
\begin{align}
\hat{\tau}_{k,n}&=\frac{2d_{k,n}}{c}+z_{\tau_{k,n}}, \\
\hat{\mu}_{k,n}&=-\frac{2(v_{k,n}\cos{\theta_{k,n}}-v_{n}^u\cos{(\theta_{k,n}+\theta^{u}_n)})f_c}{c}+z_{\mu_{k,n}}.
\end{align}
Here, $v_{k,n}$ represents the velocity of the object $k$, $v_{n}^u$ denotes the UAV's velocity, $\theta^{u}_n$ is the UAV's orientation angle relative to the horizontal plane, and $f_c$ is the carrier frequency. It is important to note that while the UAV measures radial velocity $v^R_{k,n}$ via Doppler shift, the actual velocity $v_{k,n}$ is inferred using the formula $v_{k,n} = \left(v_{n}^u\cos{(\theta_{k,n} + \theta^{u}_n)} - v^R_{k,n}\right)/\cos{\theta_{k,n}}$.

For angle measurement between the aerial base station and object $k$, methodologies such as the Capon method or the generalized likelihood ratio test (GLRT) \cite{4350230} are employed, leading to
\begin{equation}\label{estitheta}
\hat{\theta}_{k,n}=\theta_{k,n}+z_{\theta_{k,n}}.
\end{equation}
The terms $z_{\tau_{k,n}}$, $z_{\mu_{k,n}}$, and $z_{\theta_{k,n}}$ signify Gaussian measurement noises with zero mean and variances $\sigma^2_{\tau_{k,n}}$, $\sigma^2_{\mu_{k,n}}$, and $\sigma^2_{\theta_{k,n}}$, respectively. To address the complexities in determining variances for $\hat{\tau}_{k,n}$, $\hat{\mu}_{k,n}$, and $\hat{\theta}_{k,n}$, the CRB is implemented as a sensing metric for its unbiased estimation capabilities and minimum mean square error (MMSE) lower bound \cite{9705498}:
\begin{align}
\sigma^2_{\tau_{k,n}}&=CRB_\tau=\frac{1}{SNR_{k,n}N_t N_r\kappa^2},\\
\sigma^2_{\mu_{k,n}}&=CRB_\mu=\frac{1}{SNR_{k,n}N_t N_r\iota ^2},\\
\sigma^2_{\theta_{k,n}}&=CRB_\theta=\frac{1}{SNR_{k,n}N_t N_r\xi^2}.
\end{align}
Following \cite{9705498}, the squared effective bandwidth is denoted as $\kappa^2$, squared effective pulse length as $\iota^2$, and the root mean square aperture width of the beampattern is defined by
\begin{equation}
\xi^2=\frac{\pi^2 d^2_k \cos^2{\theta_{k,n}}(N_t^2-1)}{3\lambda^2}.
\end{equation}

The overall measurement model is concisely summarized as:
\begin{equation}\label{mea_h}
\left\{
\begin{aligned}
     &\hat{\theta}_{k,n}=\theta_{k,n}+z_{\theta_{k,n}},  \\
     &\hat{\tau}_{k,n}=\frac{2d_{k,n}}{c}+z_{\tau_{k,n}},\\
     &\hat{\mu}_{k,n}=-\frac{2(v_{k,n}\cos{\theta_{k,n}}-v_{n}^u\cos{(\theta_{k,n}+\theta^{u}_n)})f_c}{c}+z_{\mu_{k,n}}.
\end{aligned}
\right.
\end{equation}

\subsection{Evolution Model}

\begin{figure}[t] 
\centerline{\includegraphics[width=0.45\textwidth]{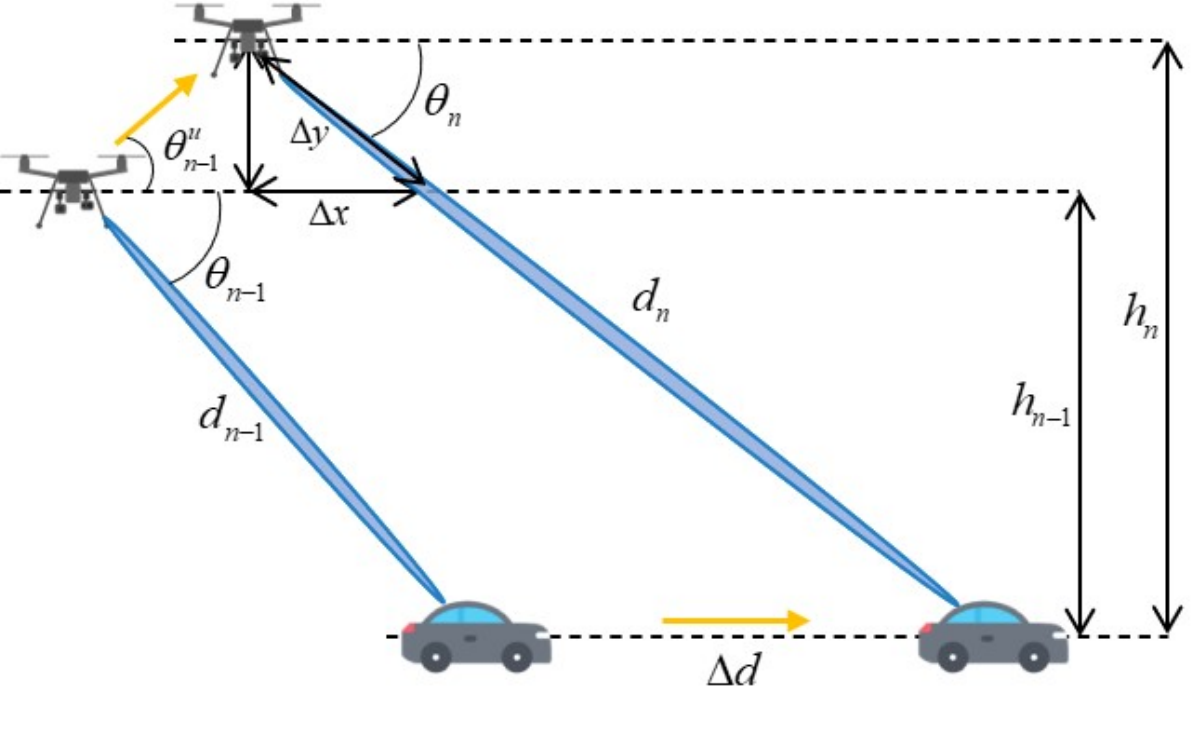}}
\caption{Evolution model.}
\label{fig.evolution}
\vspace{-0.6 cm}
\end{figure}

The evolution model is pivotal for tracking alterations in the angle, distance, and velocity of $K$ objects. For simplification, we consider a scenario where the UAV is dedicated to a single object, hence we omit the subscript $k$. The variables $\theta_n$, $d_n$, and $v_n$ represent the object's angle, distance, and velocity relative to the UAV at the $n$-th time slot, respectively.

The geometric relationship of distance is described as
\begin{align}
    d_n^2 = (h_{n-1}+ \text{\uppercase\expandafter{\romannumeral1}})^2 + (d_{n-1}\cos{\theta_{n-1}}+\text{\uppercase\expandafter{\romannumeral2}})^2,
\end{align}
where $h$ is the UAV's altitude, and
\begin{align}
    \text{\uppercase\expandafter{\romannumeral1}} &=v_{n-1}^u\sin{\theta_{n-1}^u\Delta T},\\
    \text{\uppercase\expandafter{\romannumeral2}} &=v_{n-1}\Delta T-v_{n-1}^u\cos{\theta_{n-1}^u\Delta T}.
\end{align}
Subsequent simplification yields
\begin{align}\label{evo_d}
    d_n^2 = d_{n-1}^2+ {\text{\uppercase\expandafter{\romannumeral1}}}^2 +{\text{\uppercase\expandafter{\romannumeral2}}}^2+ 2\text{\uppercase\expandafter{\romannumeral1}}  h_{n-1}+2\text{\uppercase\expandafter{\romannumeral2}} d_{n-1}\cos{\theta_{n-1}}.
\end{align}
Considering the nonlinear nature of Eq. (\ref{evo_d}), we employ an approximation
\begin{align}\label{approx_d}
    &d_n-d_{n-1} \nonumber \\
    &=\frac{{\text{\uppercase\expandafter{\romannumeral1}}}^2 +{\text{\uppercase\expandafter{\romannumeral2}}}^2+ 2\text{\uppercase\expandafter{\romannumeral1}}  h_{n-1}
    +2\text{\uppercase\expandafter{\romannumeral2}} d_{n-1}\cos{\theta_{n-1}}}{d_n+d_{n-1}} \nonumber \\
    &\approx \frac{{\text{\uppercase\expandafter{\romannumeral1}}}^2+{\text{\uppercase\expandafter{\romannumeral2}}}^2}{2d_{n-1}}+\frac{2\text{\uppercase\expandafter{\romannumeral1}} d_{n-1}\sin{\theta_{n-1}}+2\text{\uppercase\expandafter{\romannumeral2}} d_{n-1}\cos{\theta_{n-1}}}{2d_{n-1}}\nonumber \\
    &=\frac{{\text{\uppercase\expandafter{\romannumeral1}}}^2+{\text{\uppercase\expandafter{\romannumeral2}}}^2}{2d_{n-1}} + \text{\uppercase\expandafter{\romannumeral1}}\sin{\theta_{n-1}}+\text{\uppercase\expandafter{\romannumeral2}}\cos{\theta_{n-1}} \nonumber \\
    &\approx\text{\uppercase\expandafter{\romannumeral1}}\sin{\theta_{n-1}}+\text{\uppercase\expandafter{\romannumeral2}}\cos{\theta_{n-1}}.
\end{align}
where the first approximation is based on the fact that the distance of the object is very short in slot $\Delta T$, and we ignore $\frac{{\text{\uppercase\expandafter{\romannumeral1}}}^2+{\text{\uppercase\expandafter{\romannumeral2}}}^2}{2d_{n-1}}$ in the second approximation due to the square term and denominator. Accordingly, we have
\begin{align}
    d_n = d_{n-1} + \text{\uppercase\expandafter{\romannumeral1}}\sin{\theta_{n-1}}+\text{\uppercase\expandafter{\romannumeral2}}\cos{\theta_{n-1}}.
\end{align}
Focusing on the angle geometry, we derive
\begin{equation}
\left\{
\begin{aligned}
     &\frac{\Delta y}{d_n}=\frac{\Delta x}{d_{n-1}\cos{\theta_{n-1}}+\text{\uppercase\expandafter{\romannumeral2}}}=\frac{\text{\uppercase\expandafter{\romannumeral1}}}{h_n},  \\
     &\frac{\text{\uppercase\expandafter{\romannumeral2}}-\Delta x}{\sin{\Delta \theta}}=\frac{d_n-\Delta y}{\sin{\theta_{n-1}}},\\
     &h_n=h_{n-1}+\text{\uppercase\expandafter{\romannumeral1}},
\end{aligned}
\right.
\end{equation}
where $\Delta \theta = \theta_n - \theta_{n-1}$. The definitions of $\Delta x$ and $\Delta y$ are provided in Figure \ref{fig.evolution}, and they additionally adhere to the relationship $\cos{\theta_n} = \frac{\Delta x}{\Delta y}$. It's important to note that $\Delta \theta$ is of small magnitude, leading to the introduction of an approximation as follows
\begin{align}\label{approx_a}
    \Delta \theta \approx \sin{\Delta \theta}&=\frac{(\text{\uppercase\expandafter{\romannumeral2}}-\Delta x)\sin{\theta_{n-1}}}{d_n-\Delta y}\nonumber \\
    &=\frac{(\text{\uppercase\expandafter{\romannumeral2}}-\frac{(d_{n-1}\cos{\theta_{n-1}}+\text{\uppercase\expandafter{\romannumeral2}})\text{\uppercase\expandafter{\romannumeral1}}}{h_n})\sin{\theta_{n-1}}}{d_n-\frac{\text{\uppercase\expandafter{\romannumeral1}}d_n}{h_n}}\nonumber \\
    &\approx \frac{\text{\uppercase\expandafter{\romannumeral2}}\sin{\theta_{n-1}}-\text{\uppercase\expandafter{\romannumeral1}}\cos{\theta_{n-1}}}{d_{n-1}\text{\uppercase\expandafter{\romannumeral1}}\sin{\theta_{n-1}}+\text{\uppercase\expandafter{\romannumeral2}}\cos{\theta_{n-1}}}.
\end{align}
Therefore,
\begin{align}
    \theta_n=\theta_{n-1}-\frac{\text{\uppercase\expandafter{\romannumeral2}}\sin{\theta_{n-1}}-\text{\uppercase\expandafter{\romannumeral1}}\cos{\theta_{n-1}}}{d_{n-1}\text{\uppercase\expandafter{\romannumeral1}}\sin{\theta_{n-1}}+\text{\uppercase\expandafter{\romannumeral2}}\cos{\theta_{n-1}}}.
\end{align}

Under the assumption of near-constant object velocity, the model simplifies to
\begin{align}
    v_n \approx v_{n-1}.
\end{align}
The final evolution model is then presented as:
\begin{equation}\label{evo_g}
\left\{
\begin{aligned}
     &\theta_n=\theta_{n-1}-\frac{\text{\uppercase\expandafter{\romannumeral2}}\sin{\theta_{n-1}}-\text{\uppercase\expandafter{\romannumeral1}}\cos{\theta_{n-1}}}{d_{n-1}\text{\uppercase\expandafter{\romannumeral1}}\sin{\theta_{n-1}}+\text{\uppercase\expandafter{\romannumeral2}}\cos{\theta_{n-1}}}+\eta_\theta,  \\
     &d_n = d_{n-1} + \text{\uppercase\expandafter{\romannumeral1}}\sin{\theta_{n-1}}+\text{\uppercase\expandafter{\romannumeral2}}\cos{\theta_{n-1}}+\eta_d,\\
     &v_n=v_{n-1}+\eta_v,
\end{aligned}
\right.
\end{equation}
where $\eta_\theta$, $\eta_d$, and $\eta_v$ are the corresponding noises, assumed to be zero-mean Gaussian distributed with variances $\sigma_1^2$, $\sigma_2^2$, and $\sigma_3^2$, respectively. These noises stem from approximation and other systematic errors, excluding thermal noise.
\begin{figure}[h] 
\centerline{\includegraphics[width=0.45\textwidth]{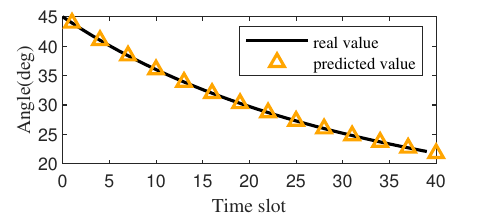}}
\centerline{\includegraphics[width=0.45\textwidth]{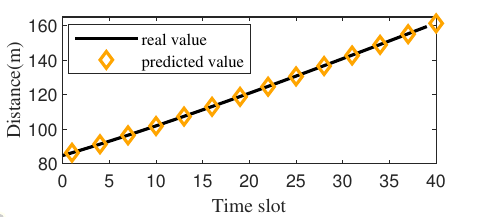}}
\caption{Validation of approximation in (\ref{approx_a}) and (\ref{approx_d}).}
\label{fig.approx}
\vspace{-0.3 cm}
\end{figure}

To evaluate the efficacy of the approximation techniques utilized for $\theta_n$ and $d_n$, we conducted a comparative analysis between actual and approximated values of angle and distance over $40$ time slots. This is illustrated in Fig. \ref{fig.approx}. The parameters selected for the UAV in this simulation were: a cruising speed of $15$m/s, an operational altitude of $60$m, and an orientation of $180^\circ$. The object being tracked had a velocity of $30$m/s and an initial angle of $45^\circ$. Each time slot in this scenario was defined as $\Delta T = 50$ms. The results indicate that despite a relatively large value of $\Delta T$, the discrepancy between real and predicted values remains minimal, suggesting that the associated variance is likely insignificant.

\section{Problem Formulation}
This section delineates the formulation of the beamforming problem in the context of sensing-assisted communication, aiming to optimize the communication rate for multiple objects. In each time slot, the UAV transmits ISAC signals to $K$ objects via $K$ distinct beams. The reflected signals are pivotal for both beam tracking, prediction, and optimization.

We adopt a temporal-assisted beamforming optimization strategy, spanning multiple time slots. Initially, during the omnidirectional phase, beams are unaltered as channel information is unavailable. For subsequent slots, beam optimization leverages sensing data from current radar echoes, with the objective of enhancing communication rates in forthcoming slots. This iterative process gradually improves sensing precision, thereby refining channel predictions and boosting the overall communication rate.

The predicted channel from the UAV to object $k$ in time slot $n$ is formulated as
\begin{equation}
\bar{\mathbf{h}}_{k,n}=\frac{\alpha_0}{\bar{d}_{k,n}}e^{j\frac{2\pi}{\lambda}d_{k,n}}\mathbf{a}(\bar{\theta}_{k,n}).
\end{equation}
The communication rate with object $k$ at time slot $n$ is described by
\begin{align}
\bar{R}_{k,n}&=\log_{2}{\left(1+\frac{|\bar{\mathbf{h}}_{k,n}^H \mathbf{w}_{k,n}|^2}{\sum_{i=1,i\neq k}^K |\bar{\mathbf{h}}_{k,n}^H \mathbf{w}_{i,n}|^2+\sigma_C^2}\right)}\nonumber \\
&=\log_{2}{\left(1+\frac{\text{tr}(\bar{\mathbf{h}}_{k,n} \bar{\mathbf{h}}_{k,n}^H \mathbf{W}_{k,n})}{\sum_{i=1,i\neq k}^K \text{tr}(\bar{\mathbf{h}}_{k,n} \bar{\mathbf{h}}_{k,n}^H \mathbf{W}_{i,n})+\sigma_C^2)}\right)},
\end{align}
where $\mathbf{W}_{k,n}$ is defined as $\mathbf{w}_{k,n}\mathbf{w}_{k,n}^H$. The optimization problem is thus expressed as
\begin{subequations}\label{opti1}
\begin{align}
&~~~~~~~~~~~~~~~~\underset{\{\mathbf{W}_{k,n}\}}{\text{max}}\sum\limits_{k=1}^K\bar{R}_{k,n} \label{opti1-a} \\
&~~~~~~~~~~~~~~~~\text{s.t.}~ \text{tr}\left(\sum\limits_{k=1}^K \mathbf{W}_{k,n}\right)\leq P_T, \label{opti1-b} \\
&~~~~~~~~~~~~~~~~\mathbf{W}_{k,n}\succeq 0,\mathbf{W}_{k,n}=\mathbf{W}_{k,n}^H,~\forall k,\label{opti1-c} \\
&~~~~~~~~~~~~~~~~\text{rank}(\mathbf{W}_{k,n})=1,~\forall k, \label{opti1-d} \\
&~~~~~~~~~~~~~~~~\bar{R}_{k,n}\geq \Gamma_k,~\forall k,  \label{opti1-e}\\
&|\mathbf{a}^H(\bar{\theta}_{k,n-1})\mathbf{W}_{k,n}\mathbf{a}(\bar{\theta}_{k,n-1})-\mathbf{a}^H(\theta_{k,n}^{\text{cover}})\mathbf{W}_{k,n}\mathbf{a}(\theta_{k,n}^{\text{cover}})|  \nonumber \\
&~~~\leq B_k \text{tr} (\mathbf{W}_{k,n}),~\forall k,~\forall |\bar{\theta}_{k,n-1}-\theta_{k,n}^{\text{cover}}|\leq l\sigma_{\theta_{k,n-1}} \label{opti1-f}.
\end{align}
\end{subequations}
In equation (\ref{opti1-a}), the aim is to maximize data throughput at the $n$-th slot. Equation (\ref{opti1-b}) enforces a transmission power constraint. The matrices $\mathbf{W}_{k,n}$ must meet semidefinite, Hermitian, and rank-one conditions as per equations (\ref{opti1-c}) and (\ref{opti1-d}). The parameter $\Gamma_k$ in equation (\ref{opti1-e}) is the minimum required communication rate for object $k$. The integration of sensing and communication is encapsulated in equation (38f), drawing inspiration from \cite{4350230}, which demonstrates the relationship between sensing error and beam width. In essence, by regulating the beam gain across potential angles, (\ref{opti1-f}) facilitates the creation of an optimized beampattern that closely approximates an ideal, adjustable-width radar beampattern. Here, $l$ typically is selected as $3$, following the Gaussian distribution principle. The term $\theta_{k,n}^{\text{cover}}$ in equation (\ref{opti1-f}) represents the collective angles satisfying the condition $|\hat{\theta}_{k,n-1}-\theta_{k,n}^{\text{cover}}|\leq l\sigma_{\theta_{k,n-1}}$, ensuring maximal angle coverage. The aim is to minimize $B_k$ for tighter and smoother beam transmission. As the sensing accuracy improves, $\sigma_{\theta_{k,n-1}}$ decreases, leading to a more focused beam, potentially augmenting the communication rate. The angular resolution is set to a precision of $0.1^\circ$, considering computational limitations.

\section{Proposed Solution}
This section introduces a Kalman filtering method for beam prediction and tracking. Given the nonlinear nature of both the measurement model and the state evolution model, a Linear Kalman Filter (LKF) is not directly applicable. Hence, we employ EKF techniques that facilitate local linearization of these nonlinear models. Furthermore, we present an algorithm aimed at deriving a high-quality solution for the optimization problem delineated in (\ref{opti1}).

\subsubsection{Beam Tracking and Prediction}
For local linearization of the nonlinear model, let the object motion parameters be represented as $\mathbf{e}=[\theta,d,v]^T$, and the measured signal parameters as $\mathbf{m}=[\theta,\tau,\mu]^T$. The models developed in (\ref{evo_g}) and (\ref{mea_h}) can be compactly rewritten as
\begin{equation}\label{}
\left\{
\begin{aligned}
     &\text{Evolution Model:}~\mathbf{e}_n=\mathbf{g}(\mathbf{e}_{n-1}) +  \boldsymbol{\eta}_n,  \\
     &\text{Measurement Model:}~\mathbf{m}_n=\mathbf{h}(\mathbf{e}_n) + \mathbf{z}_n,\\
\end{aligned}
\right.
\end{equation}
where $\mathbf{g}(\cdot)$ is as defined in (\ref{evo_g}), with noise vector $\boldsymbol{\eta} = [\eta_\theta,\eta_d,\eta_v]^T$ being independent of $\mathbf{g}(\mathbf{e}_{n-1})$. Likewise, $\mathbf{h}(\cdot)$, as per (\ref{mea_h}), includes measurement noise $\mathbf{z}_n = [z_{\theta_{n}},z_{\tau_{n}},z_{\mu_{n}}]^T$, independent of $\mathbf{h}(\mathbf{e}_n)$. Both $\boldsymbol{\eta}$ and $\mathbf{z}_n$ follow a zero-mean Gaussian distribution with covariance matrices:
\begin{align}
    &\mathbf{Q}_s=\text{diag}(\sigma_1^2,\sigma_2^2,\sigma_3^2),\\
    &\mathbf{Q}_m=\text{diag}(\sigma^2_{\theta_{n}},\sigma^2_{\tau_{n}},\sigma^2_{\mu_{n}}),
\end{align}
indicating the reliability of the predicted and measured values. The Jacobian matrices for $\mathbf{g}(\mathbf{e})$ and $\mathbf{h}(\mathbf{e})$ are essential for linearization. The Jacobian of $\mathbf{g}(\mathbf{e})$ is derived as (\ref{jacobian})
\begin{figure*}[tbp]
\begin{align}\label{jacobian}
     \frac{ \partial \mathbf{g} }{ \partial \mathbf{e} }     =\begin{bmatrix}
  1-\frac{(\text{I}\sin{\theta}+\text{II}\cos{\theta})d+\text{I}^2+\text{II}^2}{(d+\text{I}\sin{\theta}+\text{II}\cos{\theta})^2}& \frac{\text{II}\sin{\theta}-\text{I}\cos{\theta}}{(d+\text{I}\sin{\theta}+\text{II}\cos{\theta})^2} & -\frac{(d \sin{\theta}+\text{I}) \Delta T}{(d+\text{I}\sin{\theta}+\text{II}\cos{\theta})^2}\\
  v^u \Delta T \sin{\theta^u} \cos{\theta}& 1 & \Delta T \cos{\theta}\\
  0& 0 &1
\end{bmatrix}
\end{align}
\end{figure*}

The Jacobian matrices for $\mathbf{h}(\mathbf{e})$ is given by
\begin{align}
     &\frac{ \partial \mathbf{h} }{ \partial \mathbf{e} }\nonumber \\
     &=\begin{bmatrix}
  1& 0 & 0\\
  0& \frac{2}{c} & 0\\
  \frac{2f_c(v\sin{\theta}-v^u \sin{\theta + \theta^u})}{c}& 0 &-\frac{2f_c \cos{\theta}}{c}
\end{bmatrix}
\end{align}
The EKF technique is then deployed as follows:
\paragraph{Prediction of Motion Parameter}
\begin{align}
    \bar{\mathbf{e}}_{n \vert n-1}=\mathbf{g}(\bar{\mathbf{e}}_{n-1}).
\end{align}
\paragraph{Processing of Linearization }
\begin{align}
    \mathbf{G}_{n-1}=\frac{ \partial \mathbf{g} }{ \partial \mathbf{e} }\bigg|_{\mathbf{e}=\bar{\mathbf{e}}_{n-1}}, \mathbf{H}_{n}=\frac{ \partial \mathbf{h} }{ \partial \mathbf{e} }\bigg|_{\mathbf{e}=\bar{\mathbf{e}}_{n \vert n-1}}.
\end{align}
\paragraph{Prediction of MSE Matrix}
\begin{align}
    \mathbf{MSE}_{n\vert n-1} = \mathbf{G}_{n-1} \mathbf{MSE}_{n-1} \mathbf{G}_{n-1}^H + \mathbf{Q}_s.
\end{align}
\paragraph{Calculation of Kalman Gain}
\begin{align}
    \mathbf{KAL}_n =\mathbf{MSE}_{n\vert n-1}\mathbf{H}^H_n(\mathbf{H}_n\mathbf{MSE}_{n\vert n-1}\mathbf{H}^H_n+\mathbf{Q}_m)^{-1}.
\end{align}
\paragraph{Tracking of Motion Parameter}
\begin{align}
    \bar{\mathbf{e}}_{n}=\bar{\mathbf{e}}_{n \vert n-1} + \mathbf{KAL}_n(\mathbf{m}_n-\mathbf{h}(\bar{\mathbf{e}}_{n\vert n-1})).
\end{align}
\paragraph{Updating of MSE Matrix}
\begin{align}
    \mathbf{MSE}_{n} = (\mathbf{I}-\mathbf{KAL}_n\mathbf{H}_n)\mathbf{MSE}_{n\vert n-1}.
\end{align}

\textit{Note 1:} The initial MSE matrix $\mathbf{MSE}_{0}$ indicates the initial estimation uncertainty. If $\mathbf{MSE}_{0}$ is zero, it implies absolute confidence in the initial estimation. Typically, $\mathbf{MSE}_{0}$ is determined experimentally.

\textit{Note 2:} In the absence of prior information for the initial slot, the measured values are directly used as estimates.

\textit{Note 3:} In scenarios where moving objects are proximate, conventional range or angle resolution may be insufficient for effective tracking. To address this challenge, data association methods \cite{1271405, 1336475, 1413764}, such as the multiple hypothesis tracking algorithm, can be applied. These methods link and associate target data collected across different time slots, thereby enhancing the precision of echo differentiation.

\subsubsection{Beam Optimization}
The presence of interference in the objective function and the rank-one constraints in (\ref{opti1-d}) render the optimization problem non-convex, challenging the attainment of an optimal solution. We propose an approach to acquire a high-quality suboptimal solution to the problem (\ref{opti1}). This involves relaxing the rank-one constraint and employing SCA to optimize ISAC beamforming vectors, further enhanced with Semidefinite Relaxation (SDR) techniques.

Disregarding the rank-one constraint (\ref{opti1-d}), the reformulated problem (\ref{opti1}) becomes
\begin{align}
&\underset{\{\mathbf{W}_{k,n}\}}{\text{max}}\sum\limits_{k=1}^K\log_{2}{\left(1+\frac{\text{tr}(\bar{\mathbf{h}}_{k,n} \bar{\mathbf{h}}_{k,n}^H \mathbf{W}_{k,n})}{\sum_{i=1,i\neq k}^K \text{tr}(\bar{\mathbf{h}}_{k,n} \bar{\mathbf{h}}_{k,n}^H \mathbf{W}_{i,n})+\sigma_C^2)}\right)} \nonumber \\
\label{opti2}
&~~~~~~~~~~~~~\text{s.t.}~(\text{\ref{opti1-b}}),(\text{\ref{opti1-c}}),(\text{\ref{opti1-e}})~ \text{and} ~ (\text{\ref{opti1-f}}).
\end{align}

To address the non-concave nature of the objective function in problem (\ref{opti2}), SCA is employed to approximate it as a concave function through iterative processes. Denoting the iteration as $q \geq 0$ and the set $\{ \mathbf{W}_{k,n} \}$ as $\{ \mathbf{W}_{k,n}^{(q)} \}$, the following approximation is made
\begin{align}
\label{sca1}
\bar{R}_{k,n}&=\log_{2}{\left( {\textstyle \sum_{i=1}^K} \text{tr}(\bar{\mathbf{h}}_{k,n} \bar{\mathbf{h}}_{k,n}^H \mathbf{W}_{i,n})+\sigma_C^2 \right)}\nonumber \\
&~~~~-\log_{2}{\left( {\textstyle \sum_{i=1,i\neq k}^K} \text{tr}(\bar{\mathbf{h}}_{k,n} \bar{\mathbf{h}}_{k,n}^H \mathbf{W}_{i,n})+\sigma_C^2 \right)} \\
&\ge \log_{2}{\left( {\textstyle \sum_{i=1}^K} \text{tr}(\bar{\mathbf{h}}_{k,n} \bar{\mathbf{h}}_{k,n}^H \mathbf{W}_{i,n})+\sigma_C^2 \right)}\nonumber \\
\label{sca2}
&~~~~-\log_{2}{\left( {\textstyle \sum_{i=1,i\neq k}^K} \text{tr}(\bar{\mathbf{h}}_{k,n} \bar{\mathbf{h}}_{k,n}^H \mathbf{W}_{i,n}^{(q)})+\sigma_C^2 \right)} \nonumber \\
&~~~~-{\textstyle \sum_{i=1,i\neq k}^K} \text{tr} \left(\mathbf{X}_{k,n}^{(q)}(\mathbf{W}_{i,n}-\mathbf{W}_{i,n}^{(q)})\right) \triangleq \tilde{R}_{k,n}^{(q)},
\end{align}
where $\mathbf{X}_{k,n}^{(q)}$ is defined as
\begin{align}
\mathbf{X}_{k,n}^{(q)}=\frac{\log_{2}{(e)} \bar{\mathbf{h}}_{k,n}\bar{\mathbf{h}}_{k,n}^H}{{\textstyle \sum_{i=1,i\neq k}^K} \text{tr}(\bar{\mathbf{h}}_{k,n} \bar{\mathbf{h}}_{k,n}^H \mathbf{W}_{i,n}^{(q)})+\sigma_C^2}.
\end{align}

The first-order Taylor expansion is applied to the second term of  (\ref{sca1}), and its lower bound is taken, converting the objective function in (\ref{sca2}) to a convex form. Thus, problem (\ref{opti2}) evolves into the following iteration-dependent problem (\ref{opti3}):
\begin{align}\label{opti3}
\underset{\{\mathbf{W}_{k,n}\}}{\text{max}}&\sum\limits_{k=1}^K\tilde{R}_{k,n}^{(q)} \nonumber \nonumber \\
\text{s.t.}~&(\text{\ref{opti1-b}}),(\text{\ref{opti1-c}}),(\text{\ref{opti1-e}})~ \text{and} ~ (\text{\ref{opti1-f}}).
\end{align}

Problem (\ref{opti3}), now convex, is solvable via convex optimization tools like CVX. After multiple iterations, we derive an approximate iterative solution $\{ \mathbf{W}_{k,n}^{(q^\star)} \}$. If this solution satisfies the rank-one constraint, it serves as a solution to the original problem (\ref{opti1}). Typically, however, $\{ \mathbf{W}_{k,n}^{(q^\star)} \}$ does not conform to a rank-one matrix. Thus, IRM is utilized to construct an approximate rank-one solution. Let $\{ \mathbf{M}_{k,n}^{(0)} \}$ and $\{ \mathbf{M}_{k,n} \}$ represent $\{ \mathbf{W}_{k,n}^{(q^\star)} \}$ and $\{ \mathbf{W}_{k,n}^{(q^\star-1)} \}$, respectively. The necessary and sufficient condition for $\mathbf{M}_{k,n}^{(0)}$ to be rank-one is $r\mathbf{I}_{N_t-1} - \left( \mathbf{V}_{k,n}^{(0)} \right)^T \mathbf{M}_{k,n}^{(0)} \mathbf{V}_{k,n}^{(0)} \succeq 0$. Here, $\mathbf{V}_{k,n}^{(0)}$ consists of eigenvectors corresponding to $N_t-1$ smaller eigenvalues of $\mathbf{M}_{k,n}^{(0)}$, and $r$ is a positive number approaching 0. IRM iteratively reduces the rank of $\mathbf{M}_{k,n}^{(0)}$, where at the $p$-th iteration, IRM problem is formulated as
\begin{align}\label{opti4}
&~~~~~~~~~~~~~~~~\underset{\{\mathbf{M}_{k,n}^{(p)},r_p\}} {\text{max}}\sum\limits_{k=1}^K \bar{R}_{k,n}^{(p)} +w^p r_p\nonumber \\
&~~~~~~~~~~~~~~~~\text{s.t.}~ \text{tr}\left(\sum\limits_{k=1}^K \mathbf{M}_{k,n}^{(p)}\right)\leq P_T,  \nonumber\\
&~~~~~~~~~~~~~~~~\mathbf{M}_{k,n}^{(p)}\succeq 0,\mathbf{M}_{k,n}^{(p)}=\left( \mathbf{M}_{k,n}^{(p)}\right)^H,~\forall k,\nonumber \\
&~~~~~~~~~~~~~~~~\bar{R}_{k,n}^{(p)}\geq \Gamma_k,~\forall k, \nonumber \\
&|\mathbf{a}^H(\bar{\theta}_{k,n-1})\mathbf{M}_{k,n}^{(p)}\mathbf{a}(\bar{\theta}_{k,n-1})-\mathbf{a}^H(\theta_{k,n}^{cover})\mathbf{M}_{k,n}^{(p)}\mathbf{a}(\theta_{k,n}^{cover})|  \nonumber \\
&~~~\leq B_k \text{tr} (\mathbf{M}_{k,n}^{(p)}),~\forall k,~\forall |\bar{\theta}_{k,n-1}-\theta_{k,n}^{cover}|\leq l\sigma_{\theta_{k,n-1}} \nonumber \\
& ~~~~~~~~~~~~~~r_p\mathbf{I}_{N_t-1}-\left( \mathbf{V}_{k,n}^{(p-1)} \right) ^T \mathbf{M}_{k,n}^{(p)} \mathbf{V}_{k,n}^{(p-1)} \succeq 0,
\end{align}
with $\bar{R}_{k,n}^{(p)}$ expressed as
\begin{align}
\bar{R}_{k,n}^{(p)}&=\log_{2}{\left( {\textstyle \sum_{i=1}^K} \text{tr}(\bar{\mathbf{h}}_{k,n} \bar{\mathbf{h}}_{k,n}^H \mathbf{M}_{i,n}^{(p)})+\sigma_C^2 \right)}  \nonumber \\
&~~~~-\log_{2}{\left( {\textstyle \sum_{i=1,i\neq k}^K} \text{tr}(\bar{\mathbf{h}}_{k,n} \bar{\mathbf{h}}_{k,n}^H \mathbf{M}_{i,n})+\sigma_C^2 \right)} \nonumber \\
&~~~~ -{\textstyle \sum_{i=1,i\neq k}^K} \text{tr} \left( \mathbf{X}_{k,n}(\mathbf{M}_{i,n}^{(p)}-\mathbf{M}_{i,n})\right),
\end{align}
and $\mathbf{X}_{k,n}$ is defined as
\begin{align}
\mathbf{X}_{k,n}=\frac{\log_{2}{(e)} \bar{\mathbf{h}}_{k,n}\bar{\mathbf{h}}_{k,n}^H}{{\textstyle \sum_{i=1,i\neq k}^K} \text{tr}(\bar{\mathbf{h}}_{k,n} \bar{\mathbf{h}}_{k,n}^H \mathbf{M}_{i,n})+\sigma_C^2}.
\end{align}
In problem (\ref{opti4}), $w^p$ is the weight coefficient of $r_p$, increasing with the number of iterations. IRM initially calculates $\mathbf{V}_{k,n}^{(0)}$ based on $\mathbf{M}_{k,n}^{(0)}$, enters the iterative problem (\ref{opti4}) to obtain $\mathbf{M}_{k,n}^{(p)}$, and subsequently computes $\mathbf{V}_{k,n}^{(p)}$. The iterations continue until $r_p$ approaches a sufficiently small value, resulting in $\{ \mathbf{M}_{k,n}^{(p)} \}$ approaching a rank-one solution.

\section{Numerical Results}

This section presents numerical results to showcase the performance of the proposed sensing-assisted communication scheme. By default, the initial flight altitude of the UAV is set to $h = 100$m, with a vertical downward orientation, and a velocity of $15$m/s. The antenna spacing is configured as $\Delta d = \lambda/2$, and both the number of transmitting and receiving antennas for the UAV are $N_t = N_r = 30$. The variances of communication and radar noise are equal, denoted as $\sigma ^2 = \sigma_C^2 = 1$\cite{9171304}. In scenarios with multiple objects, two objects are considered with initial parameters as follows: one with $\theta_0 = 75^\circ$, $v_0 = 30$m/s, and $a = -5$m/s², and the other with $\theta_0 = 135^\circ$, $v_0 = -5$m/s, and $a = 1$m/s² in the positive direction. For single-object scenarios, the initial parameters are set to $\theta_0 = 60^\circ$, $v_0 = 30$m/s, and $a = -5$m/s². Additionally, $\varepsilon = 120(1+j)$\cite{9171304}. Key simulation parameters are summarized in Table \ref{parameter}.
\begin{table}[h]
\centering
\renewcommand{\arraystretch}{1.2}
\caption{Simulation Parameters}

\begin{tabular}{c|c}
 \hline 
\textbf{parameter} & \textbf{value}  \\
\hline\hline
Number of transmit antennas $N_t$ & 30   \\
\hline
Effective bandwidth $\kappa$ (MHz)\cite{6619570}& 80    \\
\hline
Effective pulse length $\iota$ (ms) \cite{9399786} & 0.01    \\
\hline
Signal processing gain $G_m$\cite{9171304} & 10   \\
\hline
Reference communication channel coefficient $\alpha_0$ \cite{9171304}& 1    \\
\hline
Carrier frequency $f_c$ (GHz) \cite{9171304}& 30   \\
\hline
Additive white Gaussian noise $\sigma_C^2$\cite{9171304} & 1   \\
\hline
Standard deviation of angle evolution noise  $\sigma_1$ (deg) \cite{9171304}& 0.02  \\
\hline
Standard deviation of distance evolution noise $\sigma_2$ (m) \cite{9171304}& 0.2   \\
\hline
Standard deviation of velocity evolution noise $\sigma_3$ (m/s) \cite{9171304} & 0.5   \\
\hline
Time slot length $\Delta T$ (s)& 0.01   \\
\hline

\end{tabular}
\label{parameter}
\end{table}

\begin{figure}[h] 
\centerline{\includegraphics[width=0.5\textwidth]{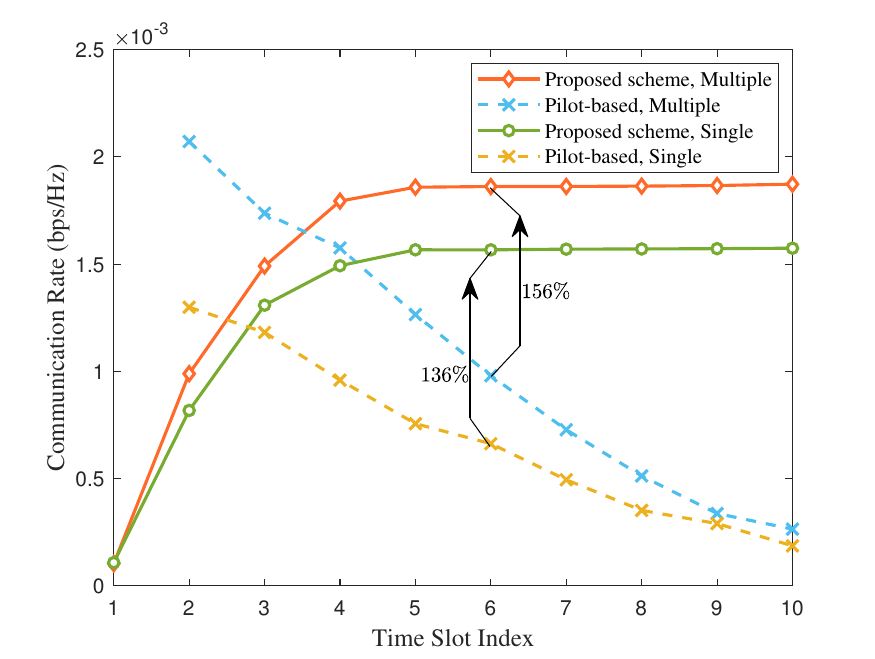}}
\caption{Comparison of communication performance between the proposed scheme and pilot-based scheme in multi-object and single-object scenarios. For UAV, $h=100$m, $v^u=15$m/s. For multiple, one with $\theta_0 = 75^\circ$, $v_0 = 30$m/s, and $a = -5$m/s², and the other with $\theta_0 = 135^\circ$, $v_0 = -3$m/s, and $a = 1$m/s². For single, $\theta_0 = 60^\circ$, $v_0 = 30$m/s, and $a = -5$m/s².
}
\label{fig.t1}
\end{figure}

The optimization of the communication rate is discussed first. Fig. \ref{fig.t1} presents a comparison between the proposed scheme and the conventional pilot-based scheme under varying numbers of objects. It is observed that the pilot scheme achieves high communication performance but loses communication opportunities in the initial time slot, and the communication rate continues to decline. In contrast, the proposed scheme exhibits a lower communication rate in the omnidirectional phase but steadily optimizes and improves its performance, and soon surpasses the pilot-based scheme. For example, in slot 6, the proposed scheme improves by $156\%$ compared to the pilot-based scheme in the presence of multi-object scenarios and $136\%$ in single-object scenarios. An important reason is that the channel based on pilot estimation has lag in the motion scene, and the pilot frequency overhead must be increased to eliminate the lag, which undoubtedly occupies the communication resource. On the other hand, the proposed scheme avoids the overhead of pilot, and the channel information can be updated at almost zero cost by using radar echoes in each time slot, thus improving the communication efficiency.

\begin{figure}[h] 
\centerline{\includegraphics[width=0.5\textwidth]{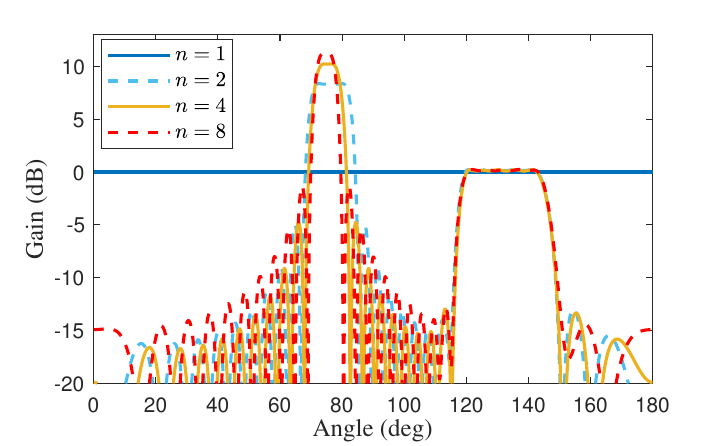}}
\caption{Optimized beampattern at different time slots in multi-object scenarios. For UAV, $h=100$m, $v^u=15$m/s. For multi-object, one with $\theta_0 = 75^\circ$, $v_0 = 30$m/s, and $a = -5$m/s², and the other with $\theta_0 = 135^\circ$, $v_0 = -3$m/s, and $a = 1$m/s².}
\label{fig.test2-1}
\vspace{-0.3 cm}
\end{figure}
\begin{figure}[h] 
\centerline{\includegraphics[width=0.5\textwidth]{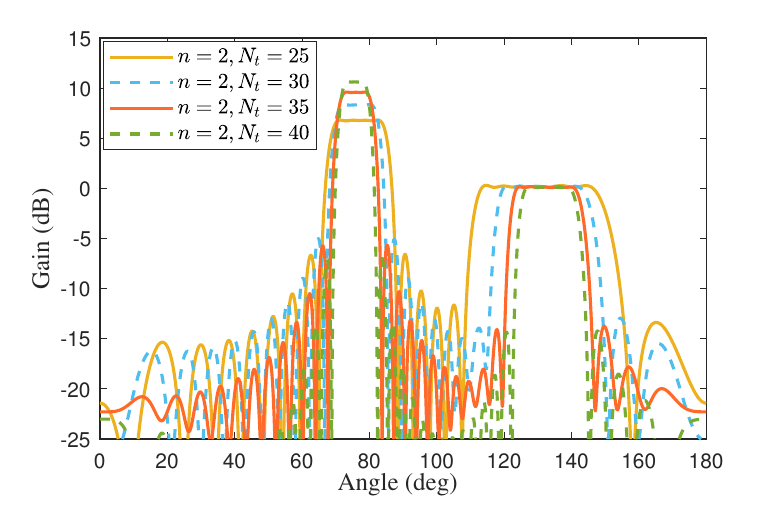}}
\caption{Optimized beampattern of different antenna numbers in multi-object scenarios when $n=2$. For UAV, $h=100$m, $v^u=15$m/s. For multi-object, one with $\theta_0 = 75^\circ$, $v_0 = 30$m/s, and $a = -5$m/s², and the other with $\theta_0 = 135^\circ$, $v_0 = -3$m/s, and $a = 1$m/s².}
\label{fig.test2-2}
\end{figure}

Fig. \ref{fig.test2-1} demonstrates the optimization results of the emission beam pattern for multiple moving objects. The UAV initially employs an omnidirectional beam pattern for the preliminary estimation of two objects. As more information is collected, the estimation becomes more accurate, resulting in a narrower beam. However, objects farther away have lower received signal-to-noise ratios, leading to wider beams and reduced power allocation. Because of the communication rate threshold, we can still guarantee the basic gain. Fig. \ref{fig.test2-2} illustrates the impact of different antenna numbers on beam optimization, showing a significant improvement in array gain with an increase in the number of transmitting and receiving antennas.

\begin{figure}[h] 
\centerline{\includegraphics[width=0.5\textwidth]{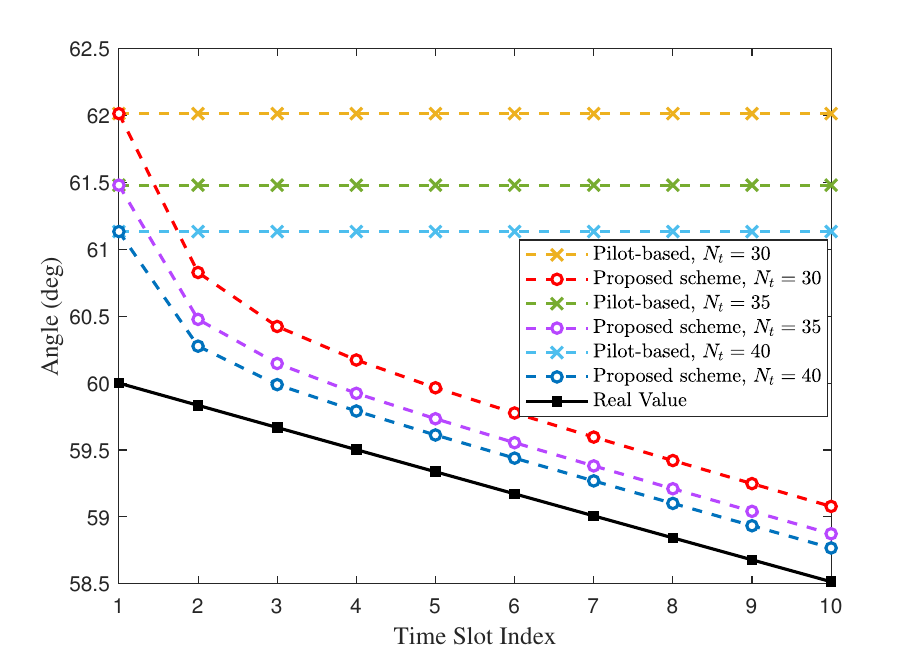}}
\caption{Comparison of angle tracking performance (Here we show it in terms of angle estimates $\bar{\theta}$.) of the proposed scheme with pilot-based scheme under different antennas. For UAV, $h=100$m, $v^u=15$m/s. For single-object, $\theta_0 = 60^\circ$, $v_0 = 30$m/s, and $a = -5$m/s².}
\label{fig.t3}
\vspace{-0.3 cm}
\end{figure}
\begin{figure}[h] 
\centerline{\includegraphics[width=0.5\textwidth]{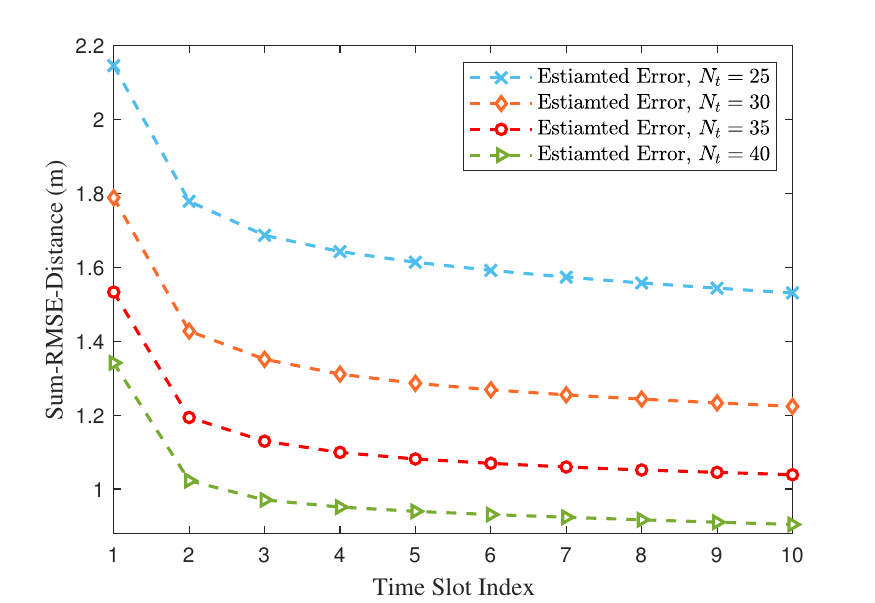}}
\caption{The tracking performance of different numbers of antennas. For UAV, $h=100$m, $v^u=15$m/s. For multi-object, one with $\theta_0 = 75^\circ$, $v_0 = 30$m/s, and $a = -5$m/s², and the other with $\theta_0 = 135^\circ$, $v_0 = -3$m/s, and $a = 1$m/s².}
\label{fig.t4}
\vspace{-0.3 cm}
\end{figure}

Angle tracking performance for single-object scenarios is shown in Fig. \ref{fig.t3}. A comparison between the proposed scheme and the pilot-based scheme reveals that it initially produces the same estimates as the proposed scheme because the matching filter gain is the same, but over time, its accuracy diminishes. The proposed scheme, on the other hand, maintains accuracy and outperforms the pilot-based scheme, with the improvement becoming more significant as the number of antennas increases.

Fig. \ref{fig.t4} demonstrates distance tracking performance using the RMSE for multiple objects tracked by the UAV. The tracking accuracy improves over time, and the addition of the EKF enables stable tracking even at high speeds. Comparing the tracking effects of different antenna numbers shows that a sufficient number of antennas ensures tracking accuracy.

\begin{figure}[h] 
\centerline{\includegraphics[width=0.5\textwidth]{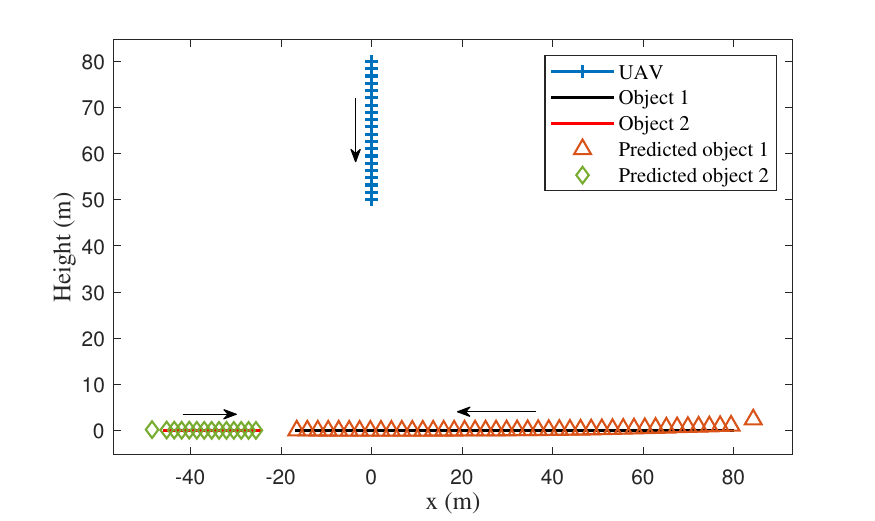}}
\caption{Dynamics of object trajectories and UAV-based predictive Analysis in multi-object environments. Note: Arrows represent the direction of movement. For UAV, $h=80$m, $v^u=8$m/s. For multi-object, one with $\theta_0 = 45^\circ$, $v_0 = -30$m/s, and $a = 0$m/s², and the other with $\theta_0 = 120^\circ$, $v_0 = 7$m/s, and $a = 0$m/s².}
\label{fig.t5}
\vspace{-0.3 cm}
\end{figure}
\begin{figure}[h] 
\centerline{\includegraphics[width=0.5\textwidth]{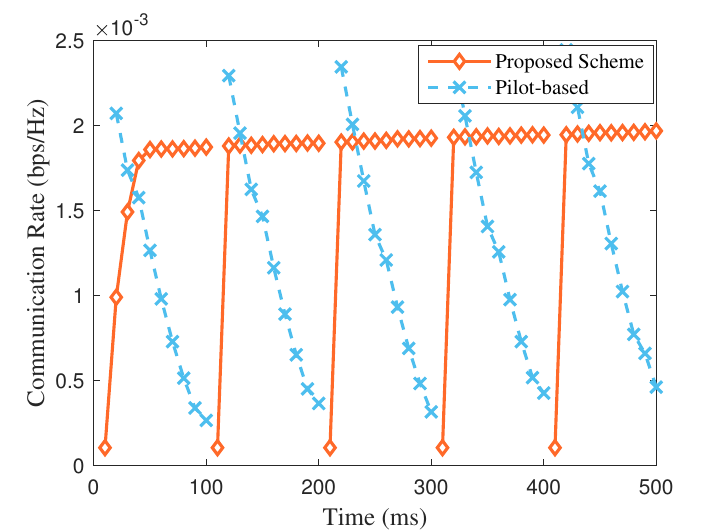}}
\caption{Communication performance for long-term tracking of multiple objects. For UAV, $h=100$m, $v^u=15$m/s. For multi-object, one with $\theta_0 = 75^\circ$, $v_0 = 30$m/s, and $a = -5$m/s², and the other with $\theta_0 = 135^\circ$, $v_0 = -3$m/s, and $a = 1$m/s².}
\label{fig.t6}
\end{figure}
Further adjustments to parameters were conducted to visualize the predictive capabilities of the UAV regarding object positions, as depicted in Fig. \ref{fig.t5}. Specifically, the initial height of the UAV was modified to $80$m, and the speed was adjusted to $8$m/s. Alterations were made to the initial angles of the two objects, setting them to $45^\circ$ and $120^\circ$, while the initial velocities were changed to $30$m/s and $7$m/s, with no acceleration setting. The position prediction effect of the EKF-based scheme, when compared to real motion trajectories, exhibited remarkable proximity to actual values after a brief optimization period. This consistency aligns with the earlier analysis and underscores the high tracking and prediction capabilities of the EKF.

\begin{figure}[h] 
\centerline{\includegraphics[width=0.5\textwidth]{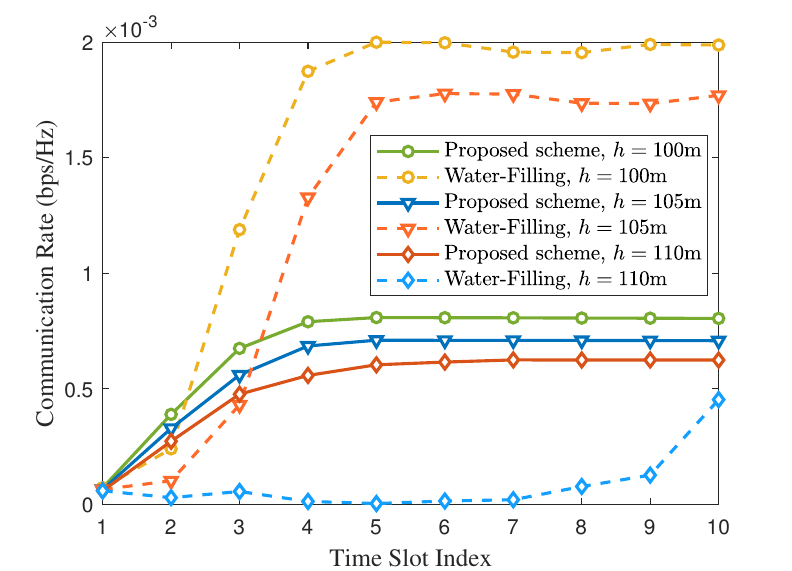}}
\caption{Comparison of communication rates between the proposed scheme and the Water-Filling scheme at different distances. For UAV, $v^u=15$m/s. For single-object, $\theta_0 = 45^\circ$, $v_0 = 30$m/s, and $a = -5$m/s².}
\label{fig.t7-1}
\end{figure}

\begin{figure}[h] 
\centerline{\includegraphics[width=0.5\textwidth]{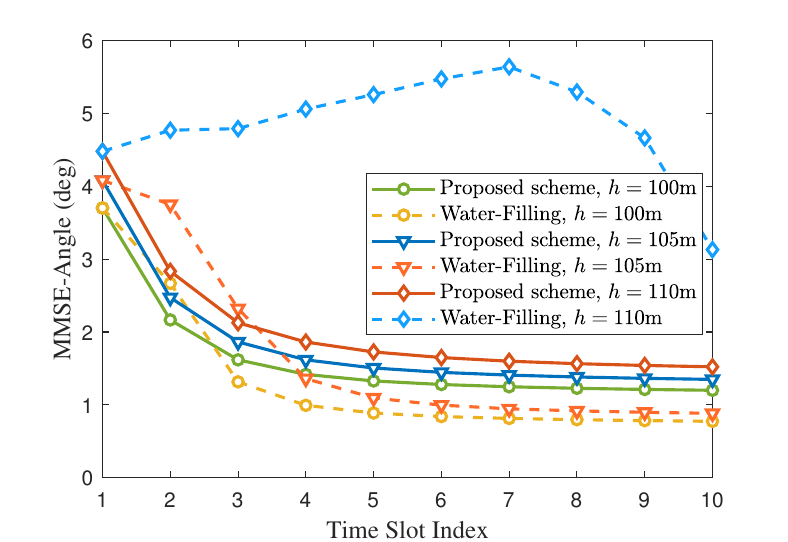}}
\caption{Comparison of angle tracking effects between the proposed scheme and Water-Filling scheme at different distances. For UAV, $v^u=15$m/s. For single-object, $\theta_0 = 45^\circ$, $v_0 = 30$m/s, and $a = -5$m/s².}
\label{fig.t7-2}
\end{figure}

Returning to the same preset parameters as in Fig. \ref{fig.t1}, we extended the time scale to $500$ms to obtain Fig. \ref{fig.t6}. The UAV periodically transmits the omnidirectional beam to update the loaded objects in the range, at the cost of intermittently losing the sensing accuracy and communication efficiency of one time slot. Fortunately, the tracking effect of EKF ensures that the communication rate quickly returns to a high level after the update. As the UAV gets closer to the objects, it can be observed that the communication rate of the Pilot-based shows a slight upward trend, and the gap between the throughput of each frame of the pilot scheme and that of the proposed scheme is also decreasing. A reasonable explanation is that the measurement error has a greater impact on the ratio of the pilot scheme. To illustrate the impact of measurement errors on the proposed scheme, more detailed comparisons are made in Figs. \ref{fig.t7-1} and \ref{fig.t7-2}.



In the following, we conducted a rough comparison between the proposed scheme and the Water-Filling scheme. In this simulation, a single object was positioned at $45^\circ$, with other conditions held constant. Fig. \ref{fig.t7-1} illustrates that as the UAV's height increased from $100$m to $110$m, the communication rates of both schemes decreased. Below $105$m, the Water-Filling scheme outperformed the proposed scheme. However, above $110$m, the Water-Filling scheme exhibited a sharp and unpredictable drop in performance, while the proposed scheme maintained stable optimization. A similar trend is confirmed by the comparison of angle tracking performance in Fig. \ref{fig.t7-2}. The Water-Filling scheme optimizes power allocation by making the beamforming vector equal to the corresponding steering vector. When an adequate number of antennas are present, the beamforming results in a ``needle-shaped" beam pattern. Compared to the proposed scheme, the ``needle-shaped" pattern exhibited clear advantages and disadvantages. When the predicted angle is sufficiently accurate, it aligns with the object, resulting in significant communication gain. However, if the predicted angle is less accurate, the main lobe becomes too narrow, causing a sharp drop in gain and severe degradation in communication performance, as evidenced in the $110$m scenario. In contrast, the proposed scheme sacrifices some gain to ensure broader beam coverage of the object, even when the predicted angle is less accurate, thereby guaranteeing a certain communication rate.
\begin{figure}[h] 
\centerline{\includegraphics[width=0.5\textwidth]{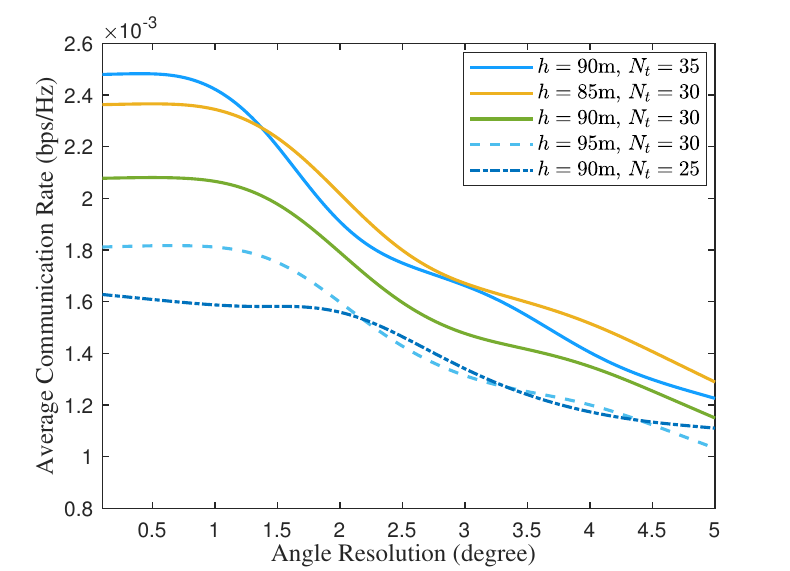}}
\caption{Average communication rate under different angular resolution. For UAV, $v^u=15$m/s. For single-object, $\theta_0 = 60^\circ$, $v_0 = 30$m/s, and $a = -5$m/s ².}
\label{fig.aa_acr}
\end{figure}

In Fig. \ref{fig.aa_acr}, we analyze the variation in the average communication rate per frame with different angular resolutions, showing that higher values of angular resolution worsen the average communication rate. With the number of antennas set at $30$ and the UAV's altitude increasing from 85m to 95m, there is no substantial decline in communication performance from $0.1^\circ$ to $1^\circ$, demonstrating minimal beamwidth adjustments and reduced power diffusion at these resolutions. At a constant altitude of $90$m, as the number of antennas increases from $25$ to $35$, the decrease in communication rate becomes pronounced earlier, starting at $1^\circ$ for $N_t=35$ and at $2^\circ$ for $N_t=25$. This trend confirms that the number of antennas is critical for precise beamwidth control, with more antennas enabling tighter focus of power on specific angles, thereby enhancing sensitivity to angular resolution variations and resulting in more pronounced changes in beamwidth with angular adjustments.

\section{Conclusion}
In this paper, we introduced a tracking and prediction scheme for the ISAC system mounted on a UAV by devising a novel frame structure. In scenarios with no prior information, the UAV initially transmits an omnidirectional beam to sensing channel to gather information. Subsequently, a directional beam is designed to enhance communication performance while obtaining more precise information for subsequent optimization. To achieve greater accuracy, we proposed an EKF method based on a newly established object motion model. The predicted values were incorporated into the optimization problem, and the SCA method was employed to transform the problem into a convex one. Subsequently, the IRM algorithm was applied to obtain high-quality solutions. Numerical results indicated that the proposed scheme outperforms the pilot-based scheme in terms of communication rate, whether in single-object or multi-object scenarios. Furthermore, a comparison with the Water-Filling scheme revealed superior optimization results for the proposed scheme, particularly at greater distances. In the future, we will focus on the research of UAV serving ground objects in 3D scenario.

\bibliographystyle{ieeetr}

\end{document}